\documentclass[12pt,leqno]{amsart}
\input psfig.sty
\newcommand\tr{{\mathrm{tr}}}
\newcommand\ad{{\mathrm{ad}}}
\newcommand\inn{{\mathrm{in}}}
\newcommand\Hom{{\mathrm{Hom}}}

\def\be{\begin{equation}}
\def\ee{\end{equation}}
\def\ba{\begin{array}}
\def\ea{\end{array}}
\def\bea{\begin{eqnarray}}
\def\eea{\end{eqnarray}}

\newcommand\cA{{\mathcal A}}
\newcommand\cB{{\mathcal B}}
\newcommand\cC{{\mathcal C}}
\newcommand\cD{{\mathcal D}}
\newcommand\cE{{\mathcal E}}

\newcommand\cG{{\mathcal G}}

\newcommand\cR{{\mathcal R}}

\newcommand\cU{{\mathcal U}}

\newcommand\BC{{\mathbb{C}}}

\newcommand\BQ{{\mathbb{Q}}}
\newcommand\BR{{\mathbb{R}}}

\newcommand\BZ{{\mathbb{Z}}}

\newcommand\fg{{\mathfrak{g}}}
\newcommand\fh{{\mathfrak{h}}}

\newcommand\fB{{\mathfrak{B}}}

\newcommand\fI{{\mathfrak{I}}}

\newcommand\fU{{\mathfrak{U}}}

\newtheorem{thm}{{Theorem}}[section]
\newtheorem{lemma}[thm]{{Lemma}}

\newtheorem{cor}[thm]{{Corollary}}

\newtheorem{remark}[thm]{{Remark}}

\newcommand{\brmk}{\begin{remark}\em}
\newcommand{\ermk}{\end{remark}}

\newcounter{itemp}
\newcounter{lastitemp}
              {%
              \begin{list}%
                    {%
                       (\roman{itemp})%
                    }%
                      {%
                          \usecounter{itemp}%
                          \setlength{\rightmargin}{\leftmargin}%
                      }%
              }%
              {%
              \setcounter{lastitemp}{\value{itemp}}%
              \end{list}%
              }%

              {%
              \begin{list}%
                    {%
                       (\roman{itemp})%
                    }%
                      {%
                          \usecounter{itemp}%
                          \setcounter{itemp}{\value{lastitemp}}%
                          \setlength{\rightmargin}{\leftmargin}%
                      }%
              }%
              {%
              \setcounter{lastitemp}{\value{itemp}}%
              \end{list}%
              }%

\begin{document}
\title{Quantum Invariants of Periodic Links and Periodic 3-Manifolds}
\author{Qi Chen and Thang Le}

\begin{abstract}
We give criteria for framed links and 3-manifolds
to be periodic of prime order. As applications
we show that the Poincare sphere is of periodicity 2, 3, 5
only and the Brieskorn sphere
(2,3,7) is of periodicity 2, 3, 7 only.
\end{abstract}

\maketitle

\noindent Math Subject Class: 57M27; 57R56

\section{Introduction}\label{intro}
For any complex simple Lie algebra $\fg$ one can define a framed
link invariant $J_L^\fg$ which is the Jones polynomial when $\fg = sl_2$.
For the definition of $J_L^\fg$, see for example \cite{le1}.
Using Chern-Simons functional and the Feynman path
integral Witten \cite{w1} introduced a 3-manifold invariant $\tau^\fg$
for any complex simple Lie algebra $\fg$. 
We will follow Reshetikhin and Turaev \cite{rt} to construct 
this 3-manifold invariant
in this note, see section~\ref{mfd1}. These invariants
are called quantum invariants because they can be defined in terms of 
the representations of the quantum group $\cU_q(\fg)$.

Murasugi \cite{m1} gave a congruence relation on $J^{sl_2}_L$
if a link $L$ is $p$-periodic. A (framed) link $L$ is
{\em $p$-periodic} if the group $H = \BZ / p\BZ$ acts on $S^3$
smoothly, with fixed point set a circle, leaving $L$ invariant.
It is also assumed that $L$ contains no fixed point. A framed
link in $S^3$ is considered as embedded annuli here.
Several authors have improved Murasugi's result in various directions
\cite{t1, y1, pr5, pr6, chb2, pr4}.
A 3-manifold $M$ is $p$-periodic if $H$ acts on $M$ smoothly with fixed
point set a circle. If $M$ is oriented then the action is required to 
be orientation preserving. We only consider 3-manifolds which are
{\em oriented, connected and closed} in this note.
For $\fg = sl_2$ Chbili \cite{chb1} and Gilmer \cite{g1} gave 
independently a necessary condition similar to
Murasugi's if $M$ is $p$-periodic for a prime $p$. The drawback of their
criterion is that it involves the quotient manifold.
In \cite{mh1, mh2} H. Murakami 
showed that if $r$ is prime and $\xi$ is a primitive $r$-th root of 
unity then the 3-manifold quantum invariant 
$\tau_\xi^{sl_2}$ essentially takes value in
$\BZ[\xi]$. Masbaum and Roberts \cite{mr}
gave a simpler proof of this fact. Based on this result Gilmer, 
Kania-Bartoszynska and Przytycki gave a necessary condition for
$r$-periodicity of integral homology spheres concerning
only $\tau_\xi^{sl_2}(M)$.

We will generalize these congruence relations, for both periodic links and 
3-manifolds, to all Lie algebras. The proof of the periodic
manifold part is made possible by the integrality of quantum invariants
of 3-manifolds \cite{le2}. 

Section~\ref{lie} recalls some basics from Lie algebra. 
Section~\ref{link} deals with the link invariant. Sections~\ref{mfd1}
and \ref{mfd2} deal with the 3-manifold invariant. In section~\ref{ap}
we show that the Poincare sphere has only prime periodicity 2, 3 and 5 and
the Brieskorn sphere $\Sigma(2, 3, 7)$ has only prime periodicity 2, 3 and
7. There we also discuss the periodicity of some other Brieskorn spheres.

\section{Basics in Lie Algebra}\label{lie}

Let $\fg$ be a simple complex Lie algebra with Cartan matrix ($a_{ij}$),
$i, j = 1, \ldots l$. Fix a Cartan subalgebra $\fh$ of $\fg$ and a set 
of basis 
roots $\alpha_1,
\ldots, \alpha_l$ in its dual space $\fh^*$. One can define a
symmetric bilinear
form ($\cdot | \cdot$) on $\fh^*$ in the following way. Multiply the $i$-th row
of ($a_{ij}$) by $d_i \in \{1, 2, 3\}$ such that ($d_i a_{ij}$) is a symmetric
matrix.
Set ($\alpha_i | \alpha_j$) $= d_i a_{ij}$. This bilinear form is 
proportional to the dual
of the Killing form restricted on $\fh$. 
Let $X$ and $Y$ be the weight lattice and root
lattice of $\fg$. 
The Weyl group $W$ acts on $X$ and $Y$ naturally. The
order of the group $G = X/Y$ is det($a_{ij}$). 

Let $X_+$ be the set of dominant weights and $Y_+ = Y \cap X_+$. 
According to the general theory of Lie algebra, the
finite-dimensional simple representations of $\fg$ 
are parameterized by the dominant weights.
Let $\cU_q(\fg)$ be the quantum group
associated to $\fg$ as defined in \cite{j} with the exception
that our $q$ is equal to $q^2$ in \cite{j}. The finite-dimensional
representations of type 1 for $\cU_q(\fg)$ are also
parameterized by the dominant weights: for every 
$\lambda \in X_+$, there is a unique finite-dimensional simple
$\cU_q(\fg)$-module $\Lambda_\lambda$ of type 1 associated to it.
Let $\cC$ (resp. $\cC^R$) be the set of all finite dimensional simple
$\cU_q(\fg)$-modules of type 1 associated to the elements in $X_+$ 
(resp. $Y_+$). 
Then $\cC^R$ is a subset
of $\cC$.

If $L$ is a (framed) link in $S^3$, a coloring (resp. root-coloring) 
of $L$ is an assignment to each of its component an element in $\cC$ 
(resp. $\cC^R$). 
Denote by $\cC_L$ (resp. $\cC_L^R$) 
the set of all colorings (resp. root-colorings) 
of $L$. If $L$ is $p$-periodic
then a coloring $c\in \cC_L$ 
is said to be $p$-periodic if all link components in one
orbit are assigned the same element in $\cC$. Denote by $\cC_L^p$ 
the set of $p$-periodic colorings. Set $\cC_L^{R,p} = \cC_L^R \cap \cC_L^p$.
Notice that a $p$-periodic coloring $c$ 
induces a coloring $c'$ on the quotient link.
If a link $L$ has color $c$, then the invariant is denoted
as $J_L^\fg(c)$.

Several constants are used
frequently in this note. Let $d = \max\{d_i\}$. 
Denote by $D$ the least positive
integer such that $D (\lambda | \mu) \in \BZ$ for all weights $\lambda$
and
$\mu$. The Coxter number is $h = 1 + (\alpha_0 | \rho)$, where $\alpha_0$
is the highest short root and $\rho$ is half of the sum of positive roots.
The dual Coxter number is $h^\vee = 1 + \max_{\alpha > 0} (\alpha | \rho)/d$.
For the exact values of these constants see, for example, \cite{le1}.

\section{Quantum Invariants of Periodic Links}\label{link}
\begin{figure}[ht]
\centerline{
\psfig{figure=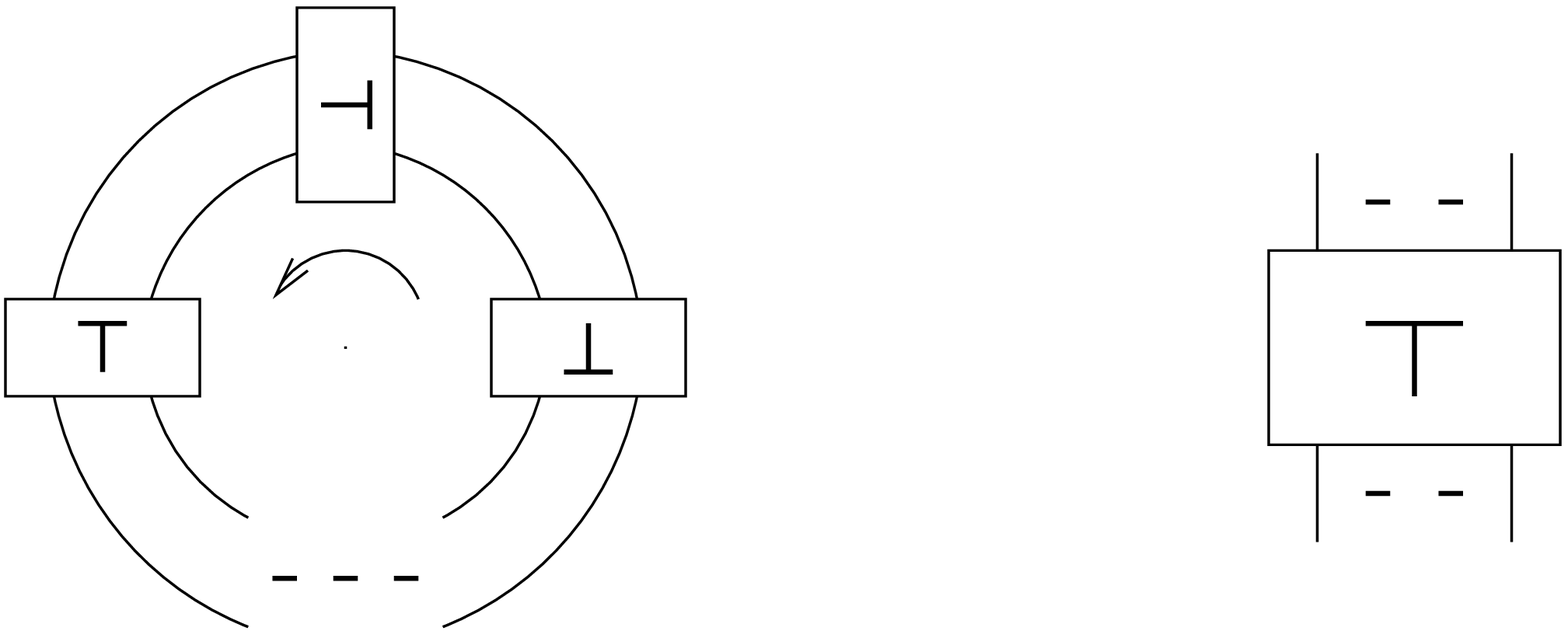,height=3cm}
}
\vskip .2cm
\caption{}\label{fig}
\end{figure}

We fix $p$ to be a prime integer throughout this paper.
If a framed link $L \subset S^3$ is $p$-periodic then $L$ 
has a diagram in $\BR^2$ with blackboard framing such that the rotation
by an angle $2\pi/p$ about
a point away from the diagram leaves the link diagram invariant
thanks to the positive answer to the Smith conjecture (\cite{mb}).
In what follows we do not distinguish a link and its diagrams if there is no
confusion.
Let $L' = L/H$ be the quotient link of $L$ with respect to the action. 
The framing on $L'$ is the blackboard framing of the quotient diagram.
There exists an ($n, n$)-tangle $T$ such that $L'$
and $L$ are the natural closure of $T$ and $T^p$ respectively.
Here $T^p$ is the tangle obtained by gluing $p$ copies of $T$ in a
natural way. See figure~\ref{fig}
where the arrow denotes the rotation generating the group $H$.

Let $U$ be the trivial knot with framing 0. The quantum dimension 
$\dim_q\Lambda_\lambda$ is defined to be $J_U^\fg(\Lambda_\lambda)$.
Let
$$
I'_{\fg,p}=\left(p, 
\ (\dim_q\Lambda_\lambda)^p-\dim_q\Lambda_\lambda, 
\ \forall \lambda \in X_+\right)
$$
be an ideal in $\BZ [q^{\pm \frac{1}{2}}]$. 
Note that $\dim_q\Lambda_\lambda$ is {\em a priori}
in $\BZ [q^{\pm \frac{1}{2D}}]$, but it is actually
an element in $\BZ [q^{\pm \frac{1}{2}}]$, see the proof of 
lemma~\ref{lm_I'inI}.
Let
$I_{\fg,p}$ be an ideal in $\BZ [q^{\pm \frac{1}{2}}]$ defined by
$$
I_{\fg,p} = \left\{ \ba{ll}
\left(p, 
\ (q^{\frac{1}{2}}+q^{-\frac{1}{2}})^p
-(q^{\frac{1}{2}}+q^{-\frac{1}{2}})\right)
= \left(p, \ (1-q^{\frac{p-1}{2}})(1-q^{\frac{p+1}{2}})\right) 
& \textrm{if $D$ is even},\\
\left(p, \ (q + q^{-1})^p-(q+q^{-1})\right) 
= \left(p, \ (1-q^{p-1})(1-q^{p+1})\right) & \textrm{if $D$ is odd}.
\ea \right.
$$
Let 
$$
I_p = \left(p, \ (q+q^{-1})^p-(q+q^{-1})\right)
= \left(p, \ (1-q^{p-1})(1-q^{p+1})\right)
$$
be an ideal in $\BZ[q^{\pm 1}]$.
One has $I_p = I_{\fg,p} \cap \BZ[q^{\pm 1}]$ when $D$ is odd.
The next lemma tells us the relationship between $I'_{\fg,p}$ 
and $I_{\fg,p}$.

\begin{lemma}\label{lm_I'inI}
$I'_{\fg,p} \subset I_{\fg,p}$.
\end{lemma}

\begin{proof}
By the strong integrality of the quantum link invariant (See \cite{le1}.), 
$\dim_q\Lambda_\lambda=J_U^\fg(\Lambda_\lambda) \in q^{(\rho|\lambda)}
\BZ[q^{\pm 1}]$.
Since $(\rho|\lambda)$ belongs to $\frac{1}{D}\BZ \cap \frac{1}{2}\BZ$,
$(\rho|\lambda)$ is an integer if $D$ is odd.
Recall that if $\bar{L}$ is the mirror image of $L$ then
$J_L^\fg$ is equal to $J_{\bar{L}}^\fg$ by substituting $q$ with $q^{-1}$. 
Since $U = \bar{U}$ as framed links one has that
$\dim_q\Lambda_\lambda=J_U^\fg(\Lambda_\lambda)$ is in 
$\BZ[q^{\frac{1}{2}}+q^{-\frac{1}{2}}]$ and is in 
$\BZ[q+q^{-1}]$ when $D$ is odd.
Recall that $p$ is prime. One has that for any $z$
$$
f^p(z) \equiv f(z) \mod (p, z^p - z) \quad\forall f(z)\in \BZ[z].
$$
\end{proof}

%Also one has the formula
%$$
%J_U(\Lambda_\lambda) = \prod_{\alpha > 0} 
%\frac{q^{(\lambda + \rho | \alpha)/2} - 
%q^{-(\lambda + \rho | \alpha)/2}}{q^{(\rho | \alpha)/2} - 
%q^{-(\rho | \alpha)/2}}.
%$$

Before we formulate the main theorem of this section let's fix
some notations first.
Suppose that a framed link 
$L$ has $m$ components with 
linking matrix $(l_{ij})$.
Let $c = (\mu_1, \ldots, \mu_m)$ be a coloring of $L$. 
Denote 
$$
f(c,L) = \sum_{1 \le i,j \le m} l_{ij} (\mu_i | \mu_j),\ \ 
\ u(c,L) = \sum_{1 \le i \le m} l_{ii} (\mu_i | 2 \rho)
$$ 
and
$$
v(c,L) = \sum_{1 \le i \le m} (\mu_i | 2 \rho).
$$
Let 
$$
\cA = \BZ[q^{1/2} + q^{-1/2}] \cdot \BZ [q^{\pm p/2}]
\oplus p \BZ [q^{\pm 1}]
$$
and 
$$\cB = \BZ[q + q^{-1}] \cdot \BZ [q^{\pm p}]
\oplus p \BZ [q^{\pm 1}]
\oplus (q^p - 1) \BZ[q^{\pm 1}].
$$
A typical element in $\cA$ is a finite summation
$(\sum a\,b)+p\,c$ where
$a \in \BZ[q^{1/2} + q^{-1/2}], b \in \BZ [q^{\pm p/2}]$ and 
$c \in \BZ[q^{\pm 1}]$. Elements in $\cB$ are similarly defined.
We introduce a normalization of $J_L^\fg$,
$$
\hat{J}_L^\fg(c) := q^{-\frac{1}{2}(f(c,L)+u(c,L)+v(c,L))} J_L^\fg(c).
$$
It is shown in \cite{le1} that 
$\hat{J}_L^\fg(c)$ does not have fractional
powers and does not depend on the framing. 
The following theorem asserts that if a framed link $L$ is $p$-periodic
with quotient link $L'$ then $\hat{J}_L^\fg$ has a special symmetry
and is closely related to $\hat{J}_{L'}^\fg$.

\begin{thm}\label{th_p_link}
Let $L$ be a $p$-periodic framed link in $S^3$ with quotient link $L'$.
Their linking matrices are $(l_{ij})$ and $(l'_{ij})$. Suppose
$L$ and $L'$ have $m$ and $m'$ components respectively.
If $p$ is a prime number, $c\in \cC_L^p$ is a $p$-periodic 
coloring of $L$, and $c'$ is the induced coloring on $L'$ then 
\begin{enumerate}
\item[(a)\ ]
$$
\hat{J}_L^\fg(c) \equiv q^{\frac{pv'-v}{2}} 
(\hat{J}_{L'}^\fg(c'))^p  ~\mod~ I'_{\fg,p} 
~\inn~  \BZ [q^{\pm \frac{1}{2}}]
$$
where $v=v(c,L)$ and $v'=v(c',L')$.
\item[(b)\ ]
If $p$ and $D$ are odd then $\frac{pv'-v}{2}$
is an integer and
$$
\hat{J}_L^\fg(c) \equiv q^{\frac{pv'-v}{2}} 
(\hat{J}_{L'}^\fg(c'))^p  ~\mod~ I'_{\fg,p}
\cap \BZ[q^{\pm 1}]  
~\inn~  \BZ[q^{\pm 1}].
$$
\item[(c)\ ]
If $(p, 2D)=1$ then
$$
q^{v^*}\hat{J}_L^\fg(c) \in \cA \cap \cB
$$
where $v^*\in\BZ$ takes the value $v/2$ if $v$ is
even and $(v+p)/2$ otherwise.
\end{enumerate}
\end{thm}

This theorem will be proved at the end of this section. 
The following observations are useful in applications.

\noindent {\em Remark.}
Because $I'_{\fg,p} \subset I_{\fg,p}$ and
$I_p = I_{\fg,p} \cap \BZ[q^{\pm 1}]$ when $D$ is odd, the statement~(b)
implies that for odd $p$ and $D$
$$
\hat{J}_L^\fg(c) \equiv q^{\frac{pv'-v}{2}} 
(\hat{J}_{L'}^\fg(c'))^p  ~\mod~ I_p
~\inn~  \BZ[q^{\pm 1}].
$$ 
Let $\xi$ be a $p$-th root of unity then
$$
q^{2v^*}\hat{J}_L^\fg(c)|_{q=\xi} \equiv
\hat{J}_L^\fg(c)|_{q=\xi^{-1}} ~\mod~ p ~\inn~  \BZ[\xi]
$$
if $(p, 2D)=1$. This follows from the statement~(c).

Several authors have proved this theorem for special Lie algebras. 
We include some of their results in corollary~\ref{special_cases} for 
comparison. First we need to recall the invariants used
in their papers.

A link invariant $P_n$ can be defined 
by the following relations:
\begin{enumerate}
\item[]
$P_n(L \sqcup L') = P_n(L)\,P_n(L')$,
\item[]
$q^{n/2}P_n(L_+) - q^{-n/2}P_n(L_-) = (q^{1/2} - q^{-1/2})P_n(L_0)$,
\end{enumerate}
where ($L_0, L_-, L_+$) is the standard skein triple.
Notice that $P_n$ is an invariant
for non-framed links.
The invariant $P_n$ can also be defined through
$\hat{J}^\fg$ as follows.
Let $L$ be a framed link whose components are colored 
by the fundamental representation of $U_q(sl_n)$, i.e.
the one corresponds to the fundamental representation of
$sl_n$. Denote this coloring by $c_n$. Let $\tilde{L}$ be any
framed link obtained from $L$ by changing the framings so that
$\tilde{L}\cdot\tilde{L}=0$, i.e. the sum of all the entries
of the linking matrix of $\tilde{L}$ is $0$. Then
$P_n(L) = J^{sl_n}_{\tilde{L}}(c_n)$. Hence
$$
\hat{J}_L^{sl_n}(c_n) = q^{(e_n-u(c_n,L)-v(c_n,L))/2}P_n(L),
$$
where $e_n = \sum_{1\leq i,j\leq m}l_{ij}\,(2\rho|\lambda_1)$.
Here $\lambda_1$ is the weight corresponding to the fundamental
representation. Let 
$[n] = \frac{q^{n/2} - q^{-n/2}}{q^{1/2} 
- q^{-1/2}}$
be the quantum integer.
The Jones polynomial can be defined as 
$$
\,V_L(t) = [2]^{-1}P_2(L)|_{\sqrt{q}=-1/\sqrt{t}}.
$$

\begin{cor}\label{special_cases}
With the same assumptions as those in theorem~\ref{th_p_link}
one has
\begin{enumerate}
\item[(a)\ ] $($\cite{chb2, pr4}$)$
$$
P_n(L) \equiv \left(P_n(L')\right)^p
\ \mod\ (p, [2]^p - [2])
\ \inn\  
\BZ[q^{\pm 1/2}] 
$$
and if $n$ is odd then
$$
P_n(L) \equiv \left(P_n(L')\right)^p \ \mod\ (p, [3]^p - [3])
\ \inn\  
\BZ[q^{\pm 1/2}].
$$
\item[(b)\ ] $($\cite{m1, y1}$)$
$$
V_L(t) \equiv (V_{L'}(t))^p \ \mod\  (p, \eta_p(t))
\ \inn \ \BZ[t^{\pm 1/2}],
$$
where $\eta_p(t) = \sum_{j=0}^{p-1}(-t)^j - t^{(p-1)/2}$.
\item[(c)\ ] $($\cite{y1}$)$
If $p$ is odd then
$$
V_L(t) - t^{2 lk} V_L(t^{-1}) \equiv 0 \ \mod\  (p,\ t^p-1)
\ \inn \ \BZ[t^{\pm 1/2}],
$$
where $lk = 1/2\sum_{1\leq i < j\leq m} l_{ij}$ i.e. $lk$ is the 
total linking number of $L$.
\end{enumerate}
\end{cor}

\noindent {\em Remark.}
The ideals $(p, [2]^p-[2])$ and $(p, [3]^p-[3])$ are used here
instead of $I_{\fg,p} (\supset I'_{\fg,p})$.
When $\fg = sl_n$ the constant $D$ is
equal to $n$.
Hence they are actually the same. Also notice that
$(t+1)\eta_p(t) \equiv q^{-p/2}([2]^p-[2])|_{\sqrt{q}=-1/\sqrt{t}}
\mod p$. One has to assume that $p$ is odd in the statement~(c)
because the Jones polynomial defined here uses 
$[2]=-\sqrt{t}-1/\sqrt{t}$ whose square is 0 mod 
$(2, t^2-1)$.
Yokota proved the statement~(c) when $L$ is a knot in
\cite{y1}.
We skip the proof of this corollary because it follows easily
from theorem~\ref{th_p_link}.

We begin the proof of theorem~\ref{th_p_link} with the next lemma
which is the $J_L^\fg$-version of theorem~\ref{th_p_link}.
The reason we use $\hat{J}_L^\fg$ instead of $J_L^\fg$ in
theorem~\ref{th_p_link} is that $\hat{J}_L^\fg$ has better
integrality.

\begin{lemma}\label{lm_p_link}
With the same assumptions as those in theorem~\ref{th_p_link} one has
\begin{enumerate}
\item[(a)\ ]
$$
J_L^\fg(c) \equiv (J^\fg_{L'}(c'))^p  ~\mod~ I''_{\fg,p} 
~\inn~  \BZ [q^{\pm \frac{1}{2D}}],
$$
where $I''_{\fg,p}=\left(p, 
~(\dim_q\Lambda_\lambda)^p-\dim_q\Lambda_\lambda, 
~\forall \lambda \in X_+\right)$ is an ideal in the ring
$\BZ [q^{\pm \frac{1}{2D}}]$.
\item[(b)\ ]
$$
J_L^\fg(c) \in \BZ[q^{\frac{1}{2}}+q^{-\frac{1}{2}}]
\cdot \BZ[q^{\pm \frac{p}{2D}}] \oplus p \BZ [q^{\pm \frac{1}{2D}}].
$$
\item[(c)\ ]
If $c \in \cC_L^{R,p}$ is a $p$-periodic root-coloring
then
$$
J_L^\fg(c) \equiv (J_{L'}^\fg(c'))^p  ~\mod~ I_p 
~\inn~  \BZ [q^{\pm 1}].
$$
\end{enumerate}
\end{lemma}

\begin{proof}
Let $T$ be the ($n, n$)-tangle mentioned at the beginning of this section.
There is a natural coloring
on $T$ by the restriction of $c$, which is still denoted by $c$.
Suppose the open ends in $T$ 
are colored by $V_1, \ldots, V_n$. Then $J_T^\fg(c)$ is a 
$\cU_q(\fg)$-module
homomorphism from $V_1\otimes\cdots\otimes V_n$ to itself, see for example
\cite{le1}.
According to \cite{le1} and \cite{lu},
$J_T^\fg(c)$ can be represented by a matrix over
$\BZ [q^{\pm \frac{1}{2D}}]$ by choosing the so called
tensor product basis for
$V_1\otimes\cdots\otimes V_n$. It is known that 
$V_1\otimes\cdots\otimes V_n$ can be decomposed over 
$\BQ(q^{\frac{1}{2D}})$ into a direct sum of
homogeneous parts, i.e. $V_1\otimes\cdots\otimes V_n = 
\bigoplus_{\lambda \in X_+} E_\lambda$ where 
$E_\lambda = \bigoplus \Lambda_\lambda$. Then $E_\lambda$ is isomorphic
to $\Lambda_\lambda\otimes N$ as $U_q(\fg)$-modules
where $N$ is a finite dimensional
vector space over 
$\BQ(q^{\frac{1}{2D}})$.
It is also known that $J_T^\fg(c)$ acts on $\Lambda_\lambda \otimes N$ 
as $Id \otimes R_\lambda$
where $R_\lambda$ is a matrix over $\BQ(q^{\frac{1}{2D}})$. 
For more details see \cite{le2}. 

If $x$ is an eigenvalue of $R_\lambda^j$ then it is also an eigenvalue
of $(J_T^\fg)^j$ for any positive integer $j$. 
Since $(J_T^\fg)^j$ has entries in 
$\BZ [q^{\pm \frac{1}{2D}}]$, $x$ must belong to $O$,
the ring of algebraic integers over
$\BZ [q^{\pm \frac{1}{2D}}]$ in some algebraic extension of
$\BQ(q^{\frac{1}{2D}})$. Hence
$\tr R_\lambda^j \in O \cap \BQ(q^{\frac{1}{2D}}) 
= \BQ[q^{\pm \frac{1}{2D}}]$.
Let $K = K_{\pm 2\rho}$, 
see \cite{le1} for the definition of $K_{\pm 2\rho}$. Then 
$\tr K|_{\Lambda_\lambda} = \dim_q\Lambda_\lambda$, 
the quantum dimension of $\Lambda_\lambda$, and
$K$ acts on $\Lambda_\lambda \otimes N$ as 
$K|_{\Lambda_\lambda} \otimes Id$. So $K J_T^{\fg}$
acts on $\Lambda_\lambda \otimes N$ as $K|_{\Lambda_\lambda} \otimes R_\lambda$.
Therefore
$$
J_{L'}^\fg(c') = \tr (K J_T^{\fg}) 
= \sum_{\lambda \in X_+} \dim_q \Lambda_\lambda ~\tr R_\lambda
$$
and
$$
J_{L}^\fg(c) = \tr (K J_{T^p}^\fg) 
= \sum_{\lambda \in X_+} \dim_q \Lambda_\lambda \tr (R^p_\lambda)
$$
in $\BZ [q^{\pm \frac{1}{2D}}]$. By lemma~\ref{lmtrace} below,
\begin{equation}\label{eq_J_L}
J_{L}^\fg(c) = \sum_{\lambda \in X_+} 
(\dim_q \Lambda_\lambda) (\tr R_\lambda)^p+p\,g,
\end{equation}
for some $g \in \BZ [q^{\pm \frac{1}{2D}}]$.
As noticed in lemma~\ref{lm_I'inI},
$\dim_q\Lambda_\lambda$ is in $\BZ [q^{\frac{1}{2}}+q^{-\frac{1}{2}}]$
and this proves the second statement.
\bea
(J_{L'}^\fg(c'))^p &=& 
(\sum_{\lambda \in X_+} \dim_q \Lambda_\lambda \tr R_\lambda)^p
\equiv \sum_{\lambda \in X_+} (\dim_q \Lambda_\lambda)^p
(\tr R_\lambda)^p \nonumber \\
&\equiv& \sum_{\lambda \in X_+} \dim_q \Lambda_\lambda (\tr R_\lambda)^p
\equiv J_L^\fg(c) ~\mod~ I''_{\fg,p}. \nonumber
\eea
This proves the first statement. 
If $c$ is a $p$-periodic root-coloring then $J_{L}^\fg(c)$ is in $\BZ[q^{\pm 1}]$ by
the strong integrality of the quantum link invariants (\cite{le1}).
Then the last statement can be proved by comparing-the-terms (see the remark
after lemma~\ref{lmtrace}).
\end{proof}

\begin{lemma}\label{lmtrace}
Let $\cD$ be an integrally closed domain 
and $F$ be its quotient field. Suppose $A$ is an
$n \times n$ matrix over $F$ whose eigenvalues are integral over $\cD$.
Then for any prime number $p, ~ \tr A^p - (\tr A)^p = p\,g$ for some
$g \in \cD$.
\end{lemma}

%The following is well known.
%\begin{fact}[Theorem 8.1, \cite{la}]\label{thms_polynomial}
%Suppose $R$ is an integral domain and $f \in R[x_1, \ldots, x_n]$ is a symmetric
%polynomial of degree $\delta$. Then there exists a polynomial $g \in R[x_1, \ldots, x_n]$
%of weight $\leq \delta$ such that $f = g(s_1, \ldots, s_n)$, where 
%$s_k$ is the k-th elementary symmetric polynomial in
%$x_1, \ldots, x_n$. The weight of a 
%monomial $x_1^{l_1} \cdots x_n^{l_n}$ is $l_1 + 2 l_2 + \cdots + n l_n$ and the
%weight of a polynomial is the maximum of the weights of the monomials occurring 
%in the polynomial.
%\end{fact}

\begin{proof}
Let $K$ be an algebraic closure of $F$. Denote by $K_\cD$ 
the integral closure of $\cD$
in $K$. Then $F \cap K_\cD = \cD$ since $\cD$ is integrally closed. 
If $x_1, \ldots, x_n$ are the eigenvalues of $A$ then $x_1^j, \ldots,
x_n^j$ are 
the eigenvalues of $A^j$ for any positive integer $j$. 
Since $x_k^j \in K_\cD$ by the assumption,
$\tr A^j = \sum_{k = 1}^n x_k^j$ belongs to $K_\cD \cap F = \cD$. Then 
$\tr\,A^p = \sum_{k = 1}^n x_k^p = (\sum_{k = 1}^n x_k)^p 
+ p\,g(s_1, \ldots, s_n)$ where $g$ is an integral coefficient
polynomial and $s_k$ is the $k$-th elementary symmetric polynomial in
$x_1, \ldots, x_n$ (Theorem 8.1, \cite{la}).
Let $e_j = \sum_{k = 1}^n x_k^j$. Then the $e_j$'s are related to the 
$s_j$'s by the Newton's formula (\cite{cohn1}),
$$
e_j - s_1 e_{j-1} + s_2 e_{j-2} + \cdots + (-1)^{j-1} s_{j-1}
e_1 + (-1)^j s_j j =0
$$
where $s_k = 0$ if $k > n$.
Because $e_j = \tr\,A^j \in \cD$ one has $s_j \in F$ 
whence $s_j \in K_\cD \cap F = \cD$.
Hence $g(s_1, \ldots, s_n) \in \cD$ as claimed.
\end{proof}

\noindent {\em Remark.}
We use a so called `comparing-the-terms' argument several places 
in this section.
We only write down the most complicated one in detail in the proof 
of theorem~\ref{th_p_link} below.

\begin{proof}[Proof of theorem~\ref{th_p_link}.]
Let $f=f(c,L)$, $u=u(c,L)$, $v=v(c,L)$,
$f'=f(c',L')$, $u'=u(c',L')$ and $v'=v(c',L')$. Notice that
$f$ and $f'$ are in $\frac{1}{D}\BZ$ and the others are integers. 
Because $f$ and $f'$ are actually summations over crossings one has $f=pf'$.
If the orbit of a link component of $L$ consists of only one component 
then its framing must be $p$ times the framing of its quotient component
in $L'$. On the other hand if an orbit consists of $p$ link components 
then they must have the same framing and color.
Therefore $u=pu'$. Then the statement~(a) follows from the statement~(a) 
of lemma~\ref{lm_p_link} and the comparing-the-terms argument.

If the orbit of a component in $L$ consists of $p$ components then
its contribution in $\frac{pv'-v}{2}$ is zero. Otherwise the contribution
is $\frac{p-1}{2}$ which is an integer for odd $p$. Now the statement~(b)
can be proved by using the comparing-the-terms argument again.
Here one must use the fact that the generators of $I'_{\fg,p}$ are in $\BZ[q^{\pm 1}]$
for odd $D$.

By comparing-the-terms and lemma~\ref{lm_p_link} we can show
that $q^{v^*}\hat{J}_L^\fg(c) \in \cA$. So we only have
to prove $q^{v^*}\hat{J}_L^\fg(c) \in \cB$. The proof is again
comparing-the-terms.
It is clear from the definition that $p|2 (v^* - v/2)$. 
Let $w = 2 (v^* - (f+u+v)/2)$ and $w' = w/p$. 
Notice that $w$ and $w'$ are in $\frac{1}{D}\BZ$.
By equation~(\ref{eq_J_L})
\bea\label{eqn_compare}
q^{v^*}\hat{J}_L^\fg & = & 
q^{w/2} J_L^\fg = q^{p w'/2} \left(\sum_{\lambda \in X_+} \dim_q \Lambda_\lambda 
(\tr R_\lambda)^p + p\,g\right) \nonumber\\
& = & \sum_{\lambda \in X_+} \dim_q \Lambda_\lambda (q^{w'/2} 
\tr R_\lambda)^p
+ p\,g'.
\eea
We split $q^{w'/2} \tr R_\lambda$ and $\dim_q \Lambda_\lambda$
into integral and fractional parts as follows.
Let $q^{w'/2} \tr R_\lambda = \alpha_\lambda + \beta_\lambda$ and
$\dim_q \Lambda_\lambda = \gamma_\lambda + \delta_\lambda$
where $\alpha_\lambda \in \BZ [q^{\pm 1}], \gamma_\lambda 
\in \BZ [q + q^{-1}],
\beta_\lambda \in \fI_D$ and $\delta_\lambda \in \fI_1 
\cap \BZ [q^{1/2} + q^{-1/2}]$
where 
$$
\fI_n = \{\sum_{finite}z_i q^i | 
z_i\in\BZ, i\in\frac{1}{2n}\BZ\backslash\BZ\}
\cup \{0\}
$$
for any positive integer $n$.
Let $\alpha_\lambda = \sum a_{\lambda, i} q^i,
\beta_\lambda = \sum b_{\lambda, j} q^j$
where $a_{\lambda, i}, b_{\lambda, j} \in \BZ$ and 
$j \in \frac{1}{2D}\BZ \backslash \BZ$.
Let $\delta_\lambda = \sum d_{\lambda, l} 
(q^l + q^{-l})$ where $l \in \frac{1}{2}\BZ\backslash\BZ$.
Then $(q^{w'/2} \tr R_\lambda)^p = A_\lambda + B_\lambda + p\,g_\lambda$ where 
$A_\lambda=\sum a_{\lambda, i}^p q^{i p}, \ 
B_\lambda = \sum b_{\lambda, j}^p q^{j p}
= \sum c_{\lambda, k} q^k$
and $g_\lambda \in \BZ [q^{\pm \frac{1}{2D}}]$.
Because $(p, 2 D) = 1$ and $j$ is not an integer, 
$k = j p$ is not an integer either. Furthermore $p|2Dk$ because
$j \in \frac{1}{2D}\BZ$. Therefore the summand in the right side
of equation~(\ref{eqn_compare}) becomes
\bea\label{eqn_summand}
\dim_q \Lambda_\lambda (q^{w'/2} \tr R_\lambda)^p & = &
(\gamma_\lambda + \delta_\lambda)
(A_\lambda + B_\lambda + p\,g_\lambda) \nonumber \\
& = & \gamma_\lambda A_\lambda + \delta_\lambda
A_\lambda + \gamma_\lambda B_\lambda + \delta_\lambda B_\lambda +
p\,g_\lambda (\gamma_\lambda + \delta_\lambda).
\eea
One has 
$\delta_\lambda B_\lambda = \sum_{k,l} c_{\lambda, k} 
d_{\lambda, l} (q^{k+l} + q^{k-l})$.
Notice $k + l \in \BZ$ if and only if 
$k - l \in \BZ$ since $2 l \in \BZ$.
Meanwhile if $k + l \in \BZ$ then $2 k = k + l + k - l \in \BZ$.  
Hence if $k+l$ is an integer then $p|2 k$ 
because $p|2 D k$
and $(p, D) = 1$. 
Therefore $k + l \equiv -(k - l) ~\mod~ p$ if $k + l \in \BZ$.
Let $\delta_\lambda B_\lambda 
= C_\lambda + D_\lambda$ where $C_\lambda \in \BZ [q^{\pm 1}]$
and $D_\lambda \in \fI_D$. 
Then $C_\lambda$ is actually an element in 
$\BZ [q + q^{-1}] \oplus 
(q^p - 1)\BZ [q^{\pm 1}]$.

Using equation~(\ref{eqn_summand}) one can expand 
equation~(\ref{eqn_compare}) as follows.
$$
q^{v^*} \hat{J}^\fg_L = 
\sum_\lambda (\gamma_\lambda A_\lambda + C_\lambda) + 
\sum_\lambda
(D_\lambda + \delta_\lambda A_\lambda + \gamma_\lambda B_\lambda) 
+ p \sum_\lambda
(G_\lambda + G'_\lambda)
$$
where $g_\lambda = G_\lambda + G'_\lambda$ and $G_\lambda \in \BZ 
[q^{\pm 1}]$, $G'_\lambda \in \fI_D$.
Notice that $\delta_\lambda A_\lambda, \gamma_\lambda B_\lambda$ 
and  $D_\lambda$ are in $\fI_D$.
Therefore $\sum_\lambda
(D_\lambda + \delta_\lambda A_\lambda + \gamma_\lambda B_\lambda + 
p ~G'_\lambda) =0$
since $q^{w/2} J^\fg_L \in \BZ [q^{\pm 1}]$ by the strong integrality 
of the quantum link invariants 
(\cite{le1}).
This completes the proof.
\end{proof}

\section{Quantum Invariants of Periodic 3-Manifolds I}\label{mfd1}

Let $M$ be a connected, oriented and closed 3-manifold. 
Recall that for any Lie algebra $\fg$ and some integer
$r>1$ one can define a 3-manifold invariant $\tau_\xi^\fg(M)$
where $\xi$ is a primitive $r$-th root of unity.
Kirby and Melvin \cite{km} proved that when the Lie algebra
is $sl_2$ and $r$ is odd $\tau_\xi^{sl_2}(M)$ is a product of
two invariants, namely 
$\tau_\xi^{sl_2}(M) = \tau_\xi^G(M)\tau_\xi^{Psl_2}(M)$.
The first factor is a weak invariant which is determined by the first
homology group and the linking form on its torsion.
The second factor is called the projective quantum invariant
which will be defined later in this section. In \cite{le2} 
the second author has generalized
this splitting to all quantum groups.
For a $p$-periodic 3-manifold $M$ one can show that the values of
$\tau_\xi^{\fg}(M)$ and $\tau_\xi^{P\fg}(M)$ are restricted.
We will only formulate criteria for the projective quantum invariant  
since the proofs 
for the non-projective ones 
are very similar. Furthermore the projective part contains `almost
all' information as noted above. Another reason is that
the projective invariants possess nice integrality (\cite{le2}).

First let's recall briefly the definition of projective quantum 
invariants of 3-manifolds.
For more details please see \cite{le2}.
Let $\fh^*_\BR$ be the
$\BR$-vector space spanned by $\alpha_1,
\ldots, \alpha_l$. Then $\fh^*_\BR \otimes \BC = \fh^*$.
Fix an integer $r > d h^\vee$. 
Let $\xi$ be
a primitive $r$-th root of unity.
The {\em fundamental alcove of level $k$} is defined as
$$
C_k = \{x \in \fh^*_\BR ~|~ (x|\alpha_i) \ge 0, (x|\alpha_0) < k, i = 1,
\ldots, l\},
$$
where $k=r-h>0$ and $\alpha_0$ is the short highest root associated 
to the basis roots we choose.
Let 
$$
P_r = \{x = c_1\alpha_1 + \cdots + c_l\alpha_l \in \fh^*_\BR ~|~ 
0 \le c_1, \ldots, c_l < r\}.
$$
Let $\bar{C}_k$ and $\bar{P}_r$ be their closures.
Recall that $\cC^R$ is the set of all finite dimensional simple
$\cU_q(\fg)$-modules of type 1 with highest weight in the root lattice. 
Let $L$ be a framed link in $S^3$ with $m$ components.
A root-coloring 
of $L$ is an assignment to each of its component an element in $\cC^R$. 
Denote by $\cC_L^R$ the set of all root-colorings 
of $L$. One can extend the colors to all weights through the
Weyl group action on $\fh^*$ (\cite{le2}).  
By $\bar{C}_k^m \cap \cC_L^R$ (resp. $P_r^m \cap Y$)
we mean
that the color on each link component is from the set $\bar{C}_k$ (resp.
$P_r$). Define
$$
F_L^{P\fg}(\xi) = \sum_{c \in \bar{C}_k^m \cap \cC_L^R} Q^\fg_L (c) 
|_{q = \xi}.
$$
where $Q^\fg_L (c) = J^\fg_L (c) 
J^\fg_{U^{(m)}}(c)$ and $U^{(m)}$ is the trivial link with
$m$ components.
Let $W$ be the Weyl group of $\fg$.
By the first symmetry principle (\cite{le2, le1}), if $r \geq d h^\vee$
then
$$
F_L^{P\fg}(\xi) = 
\left(\frac{1}{|W|}\right)^m\sum_{c \in P_r^m \cap Y} Q^\fg_L (c) 
|_{q = \xi}.
$$
Let $U_\pm$ be the unknot with $\pm 1$ framing. 
According to \cite{le2},
\begin{equation}\label{eq_U_gauss}
F_{U_+}^{P\fg}(\xi) = \frac{\gamma^{P\fg}(\xi)}{\prod_{\alpha
> 0}(1 - \xi^{(\alpha | \rho)})},
\end{equation}
where $\gamma^{P\fg}(\xi) = \sum_{\mu \in P_r \cap Y} 
\xi^{\frac{|\mu + \rho|^2 - |\rho|^2}{2}}$ is a Gauss sum. 

We need one lemma
about the Gauss sum to understand $F^{P\fg}_{U_+}$. First let's recall some 
terminologies. Let $\cG$ be a finite abelian group.
A quadratic form $z$ on $\cG$ is a function $z: \cG \to \BQ/\BZ$
satisfying the  following conditions:
\begin{enumerate}
\item[]
$z(n g) = n^2 z(g), \forall g \in \cG$ and $n \in \BZ.$
\item[]
$b_z(g, g') := z(g + g') - z(g) - z(g')$ is a symmetric bilinear form on $\cG$.
\end{enumerate}
Denote by $\ad\,b_z : \cG \to \Hom_\BZ(\cG, \BQ/\BZ)$ the adjoint representation 
of $b_z$.

\begin{lemma}[\cite{fl1}]\label{lm_gauss}
$$
\sum_{g \in \cG} e^{2 \pi i z(g)} = \sqrt{|{\mathrm{ker ~ad}}~b_z| |\cG|} 
~\epsilon,
$$
where $\epsilon$ is $0$ or an eighth root of unity.
\end{lemma}

\begin{lemma}\label{lm_denominator}
If $r > d h^\vee$ and $F_{U_\pm}^{P\fg}(\xi) \ne 0$ then 
$1/F_{U_+}^{P\fg}(\xi) \in \BZ[\xi, \frac{1}{r}, \epsilon]$ and
$\frac{F^{P\fg}_{U_+}(\xi)}{F^{P\fg}_{U_-}(\xi)} 
= \xi^{-((r + 1)^2 + 2)|\rho|^2} \omega$
where $\epsilon$ is an eighth root of unity and
$\omega$ is a fourth root of unity.
\end{lemma}

\begin{proof}
Let $\mu \in Y$ then $(\mu | \rho) \in \BZ$.
$$
|\mu + \rho|^2 - |\rho|^2 \equiv |\mu + (r+1) \rho|^2
- (r + 1)^2 |\rho|^2 ~\mod~ 2 r.
$$
Hence
$$
\gamma^{P\fg}(\xi) = \sum_{\mu \in P_r \cap Y} 
\xi^{\frac{|\mu + (r + 1) \rho|^2 - (r + 1)^2|\rho|^2}{2}}
= \xi^{-\frac{(r + 1)^2|\rho|^2}{2}} \sum_{\mu \in P_r \cap Y} 
\xi^{\frac{|\mu + (r + 1) \rho|^2}{2}}.
$$
Let $\gamma' = \sum_{\mu \in P_r \cap Y} 
\xi^{\frac{|\mu|^2}{2}}$ and
$\gamma'_\lambda = \sum_{\mu \in P_r \cap Y} 
\xi^{\frac{|\mu+\lambda|^2}{2}}$.
It's easy to see that 
$\gamma' = \gamma'_\lambda$ for any $\lambda \in Y$.
Therefore 
$\gamma^{P\fg}(\xi) = \xi^{-\frac{(r + 1)^2|\rho|^2}{2}}
\gamma' \in \BZ[\xi]$. Notice that
$P_r \cap Y$ can be identified canonically with the abelian group
$Y/rY$ with $r^l$ elements where $l$ is the dimension of $\fh$.
Let $z: P_r \cap Y \to \BQ/\BZ$ defined by
$z(\mu) \equiv \frac{|\mu|^2}{2r} ~\mod~ \BZ$. Then
$e^{2 \pi i z(\mu)} = \xi^{\frac{|\mu|^2}{2}}$.
$$
\frac{|\mu + r \alpha_i|^2}{2 r} - \frac{|\mu|^2}{2 r} =
\frac{1}{2r}\left(2 (\mu | r \alpha_i) - r^2 (\alpha_i | \alpha_i)\right)
\in \BZ
$$
where $\alpha_i$ is a basis root.
Hence
$$
z(n \mu) \equiv \frac{|n \mu|^2}{2r} = n^2 \frac{|\mu|^2}{2 r}
\equiv n^2 z(\mu) ~\mod~ \BZ
$$
and
$$
b_z(\mu, \lambda) = z(\mu + \lambda) - z(\mu) - z(\lambda) \equiv
\frac{|\mu + \lambda|^2}{2r} - \frac{|\mu|^2}{2r} - 
\frac{|\lambda|^2}{2r} = \frac{(\mu | \lambda)}{r} ~\mod~ \BZ.
$$
Hence $z$ is a quadratic form on $P_r \cap Y$.
Then $\ker \ad b_z$ is a subgroup
of $P_r \cap Y$, and in particular $|\ker \ad b_z|$ divides 
$|P_r \cap Y| = r^l$. By lemma~\ref{lm_gauss}
and equation~(\ref{eq_U_gauss}) one has
$
F_{U_+}^{P\fg}(\xi)={\sqrt{s} \epsilon}/{f(\xi)}
$
for some integer $s$, some eighth root of unity
$\epsilon$ and some $f(\xi)\in\BZ[\xi]$. Here $s$ divides 
$r^k$ for some $k$.
On the other hand $F_{U_+}^{P\fg}(\xi)$ is in $\BZ[\xi]$ by
definition so $\sqrt{s}$ belongs to $\BZ[\xi, \epsilon]$.
This proves the first statement.
For the second statement one notices that 
$F_{U_+}^{P\fg}(\xi) = \overline{F_{U_-}^{P\fg}(\xi)}$ 
where the bar stands for the complex conjugation.
\end{proof}

If $F_{U_\pm}^{P\fg}(\xi) = 0$,
we define $\tau_\xi^{P\fg}(M) = 0$ otherwise let
$$
\tau_\xi^{P\fg}(M) := 
\frac{F_L^{P\fg}(\xi)}
{(F_{U_+}^{P\fg})^{\sigma_{L+}}(F_{U_-}^{P\fg})^{\sigma_{L-}}},
$$
where $M$ is obtained by surgery along the framed link $L$ and $\sigma_{L+}$
(resp. $\sigma_{L-}$)
is the number of positive (resp. negative) eigenvalues of the 
linking matrix of $L$.
One can show that $\tau_r^{P\fg}$ is a 3-manifold invariant and is 
called the projective quantum invariant.

A framed link $L$ is called {\em strongly $p$-periodic} if $L$ is
$p$-periodic and there is no component of $L$ invariant under the
$\BZ/p\BZ$ action.

The following theorem is proved recently by Przytycki and Sokolov.

\begin{thm}[\cite{ps}]\label{thm_ps}
Let $p$ be a prime integer and $M$ be a closed, oriented and
connected 3-manifold.
Then $M$ is $p$-periodic iff there exists a strongly $p$-periodic
framed link $L \subset S^3$ such that $M$ is the result of surgery on $L$.
\end{thm}

Denote by $R_\fg^\xi$ the ring $\BZ[\xi, \frac{1}{r}, \epsilon]$
where $\epsilon$ is
the same as in lemma~\ref{lm_denominator}.
Let $I_{\fg,p}^\xi$ be the ideal
$(p, \ (\xi + \xi^{-1})^p - (\xi + \xi^{-1}))$ in 
$R_\fg^\xi$.

\begin{thm}\label{th_p_mfd}
Let $M$ be a $p$-periodic rational 
homology 3-sphere with quotient manifold $M'$
where $p$ is a prime number. 
Let $r > d h^\vee$ be an odd integer and $\xi$
be a primitive $r$-th root of unity then
$$
\tau_\xi^{P\fg}(M) \equiv (-\xi)^u (\tau_\xi^{P\fg}(M'))^p ~\mod~ 
I^\xi_{\fg,p}
~\inn~ R_\fg^\xi
$$
for some integer $u$ if $p$ is not a factor of $r |W|$. ($W$ is the Weyl group.)
\end{thm}

\noindent {\em Remark.}
Theorem~\ref{th_p_mfd} is proved for $\fg = sl_2$ 
by Chbili \cite{chb1} and
Gilmer \cite{g1}.

\begin{proof}[Proof of theorem~\ref{th_p_mfd}.]
If $I^\xi_{\fg,p} = R_\fg^\xi$
or $F_{U_+}^{P\fg} = 0$ then the theorem holds trivially.
(For example if $r$ is prime, not equal to $p$ and
not a factor of $p^2-1$ then 
$I^\xi_{\fg,p} = R_\fg^\xi$.)
Now suppose otherwise then $|W|$ is invertible in 
$R_\fg^\xi/I^\xi_{\fg,p}$. 
By lemma~\ref{lm_denominator}, $1/F_{U_\pm}^{P\fg}$ belongs to 
$R_\fg^\xi$ so  $F_{U_\pm}^{P\fg}$ is invertible in $R_\fg^\xi$.
By theorem~\ref{thm_ps}, $M$ is the result of surgery on a strongly $p$-periodic
framed link $L \subset S^3$. Let $L'$ be the quotient link in $S^3$. Then $M'$ is the 
result of surgery on $L'$. $M'$ is also a rational homology 3-sphere because 
$H_1(M) \to H_1(M')$ is a surjection under the covering projection
(actually $\pi_1(M) \to \pi_1(M')$ is a surjection).
By lemma~\ref{lm_p_link},
$$
Q_L^\fg (c) = J^\fg_L (c) 
J^\fg_{U^{(m)}}(c) \equiv (J^\fg_{L'} (c'))^p (J^\fg_{U^{(m/p)}}(c'))^p
= (Q^\fg_{L'} (c'))^p
$$
modulo $I^\xi_{\fg,p}
~\inn~ R_\fg^\xi$
for any $c \ \inn\  \cC_L^{R,p}$.
Hence we have
\begin{eqnarray}
|W|^m F_L^{P\fg}(\xi) & = & 
\sum_{c \in P_r^m \cap Y} Q^\fg_L (c) 
|_{q = \xi} \equiv 
\sum_{c \in P_r^m\cap\cR_p} Q^\fg_L (c) 
|_{q = \xi} \nonumber\\
& \equiv & \left(\sum_{c \in P_r^{m/p} \cap Y} Q^\fg_{L'} (c) 
|_{q = \xi}\right)^p \nonumber\\
 & \equiv & |W|^m \left(F_{L'}^{P\fg}(\xi)\right)^p
~\mod~ I^\xi_{\fg,p}
~\inn~ R_\fg^\xi \nonumber
\end{eqnarray}
where $\cR_p$ is the set of $p$-periodic root colorings.
Recall that we extend the set of colors to all roots not just
dominant ones.
The first equivalence follows from the fact that non-periodic
coloring must occur $p$ multiple times. One can see this fact by letting
$\BZ/p\BZ$ act on the set of colorings of $L$.
Since $|W|$ and $F^{P\fg}_{U_\pm}(\xi)$ are invertible in
$R_\fg^\xi/I^\xi_{\fg,p}$ one has
$$
\tau^{P\fg}_\xi(M) 
(F_{U_+}^{P\fg})^{\sigma_{L+}}(F_{U_-}^{P\fg})^{\sigma_{L-}} = 
F_L^{P\fg}(\xi) \equiv (F_{L'}^{P\fg}(\xi))^p 
$$
$$
= \left(\tau_\xi^{P\fg}(M') (F_{U_+}^{P\fg})^{\sigma_{L'+}}
(F_{U_-}^{P\fg})^{\sigma_{L'-}}\right)^p,
$$
and hence
$$
\tau_\xi^{P\fg}(M) \equiv (\tau^{P\fg}_\xi(M'))^p 
(F_{U_+}^{P\fg})^{p \sigma_{L'+} - \sigma_{L+}}
(F_{U_-}^{P\fg})^{p \sigma_{L'-} - \sigma_{L-}}
~\mod~ I^\xi_{\fg,p}
~\inn~ R_\fg^\xi.
$$
Let $a = p \sigma_{L'+} - \sigma_{L+}$ and $b = p \sigma_{L'-} 
- \sigma_{L-}$ then
$a + b = p (\sigma_{L'+} + \sigma_{L'-}) - (\sigma_{L+} + \sigma_{L-})
= m - m =0$ because $M$ and $M'$ are rational homology 3-spheres.
Therefore
$$
\tau^{P\fg}_\xi(M) \equiv (\tau^{P\fg}_\xi(M'))^p 
\left(\frac{F_{U_+}^{P\fg}}
{F_{U_-}^{P\fg}}\right)^a.
$$
Let $\omega$ be the fourth root 
of unity as in lemma~\ref{lm_denominator}.
Because $r$ is odd one has $\BQ(\xi) \cap \BQ(i) = \BQ$
(see \cite{ft}), hence $\omega = \pm 1$. 
$\left(\frac{F_{U_+}^{P\fg}}
{F_{U_-}^{P\fg}}\right)^a = \pm \xi^{u'} = (-\xi)^u$
since $r$ is odd and $\xi$ is an $r$-th root of unity.
\end{proof}

\section{Quantum Invariants of Periodic 3-Manifolds II}\label{mfd2}

Theorem~\ref{th_p_mfd} does not help to
determine whether a 3-manifold is
$p$-periodic. It can only tell us
that a 3-manifold is not a $p$-fold branched covering space over a particular
3-manifold branched over a circle as noted by Chbili in \cite{chb1} where it is
showed that the Poincare sphere is not a branched cover of $S^3$ of order 11 with
the branched set a circle.
Recently Gilmer, Kania-Bartoszynska and Przytycki \cite{gkp} gave a criterion
for an integral homology 3-sphere $M$ to be $r$-periodic,
which involves only the quantum $sl_2$ invariant of $M$ at $r$-th 
root of unity.
The purpose of this section is to generalize their criterion to 
quantum invariants
of all simple Lie algebras. 

The proofs in \cite{gkp} depend on the topological quantum field theory
(TQFT) behind the quantum $sl_2$ invariants of 3-manifolds.
We recall briefly some related definitions first. See \cite{tu1, g2}
for more details. 
A TQFT ($V$, $Z$) is a functor from the 3-dimensional 
cobordism category to
the category Mod($K$) of projective modules over 
a ring $K$ where $V$ and $Z$ are
maps between objects and morphisms respectively.
The objects and morphisms in the 3-dimensional 
cobordism category are surfaces 
and 3-dimensional manifolds respectively.
Turaev \cite{tu1} showed that given a modular category one can construct
a TQFT. In order to remove the anomaly one also needs to consider manifolds
with certain structure ($\fB, \fU$).
More explicitly, a structure $\fU (\Sigma)$  
on a surface $\Sigma$ is a pair ($b$, $l$) where $b$ is a finite set of 
disjoint, directed line segments colored by
objects from a modular category
and $l$ is a Lagrangian subspace of $H_1(\Sigma, \BQ)$. A structure
$\fB (M)$
on a 3-dimensional manifold $M$ is a pair ($\Omega, w$) where $\Omega$ 
is a properly embedded ribbon graph (with coupons) in $M$ colored by objects 
and morphisms from the same
modular category
and $w$ is an integer.
This cobordism category will be denoted as $\cC(\fB, \fU)$.
Gilmer \cite{g2} introduced a notion of {\em almost integral TQFT}.
Let $K$ be an integral domain which contains a Dedekind domain $\cD$. 
Then the TQFT ($V$, $Z$) is almost ($\cD$-)integral if there is an
element $\cE \in K$ such that
$\cE<M> \in \cD$ for any connected closed 3-manifold 
$M$ where $<M> = Z(M)(1)$.

\noindent {\em Remark.}
$\cE$ is required to be in $\cD$ in \cite{g2} but the proves there work
for the more general case. 

\noindent {\em Remark.}
If the modular category comes from representations of $\cU_q(\fg)$ at root
of 1 then the TQFT constructed above is not almost integral. 
To get an almost integral TQFT one has to 
modify the structures on the cobordism category.
The new cobordism category
will be denoted by $\cC(\fB^e, \fU^e)$. See \cite{chenle2}.

Let $\kappa$ and $\eta$ be one of the
square roots of ${\frac{F_{U_-}^{P\fg}}{F_{U_+}^{P\fg}}}$ and 
${F_{U_-}^{P\fg} F_{U_+}^{P\fg}}$ respectively such that
$\kappa \eta = F_{U_-}^{P\fg}$.
Set $K = \BC$ and the Dedekind domain $\cD = \BZ[\xi, \kappa]$.

\begin{lemma}\label{lm_almost}
Let $r > d h^\vee$ be a prime and not a factor of $|G| |W|$.
Let $\xi$ be a primitive $r$-th root of unity. 
Then one can define an almost integral TQFT (Z, V) 
from $\cC(\fB^e, \fU^e)$
to Mod(K) such that $\eta <M>$ is in $\cD$ for closed, connected M. 
Furthermore, if M has
0 weight and the ribbon graph in M is empty then $<M> = \tau_\xi^{P\fg}(M)$.  
\end{lemma}

\begin{proof}[Sketch of the proof]
 In \cite{le2} the second author
showed that for any simple Lie algebra $\fg$ there exist modular tensor 
categories $\cC_\xi$
using the representations of $\cU_q(\fg)$ associated to root lattice.
One can show that
Turaev's construction of anomaly-free TQFT
can be adopted to the modified 
cobordism category $\cC(\fB^e, \fU^e)$.
See \cite{chenle2} for more details.
\end{proof}

%%%%%%%%%%%%%%%%%%%%
%\begin{proof}
%Suppose that $M$ have structure ($\Omega$, $w$) and can be obtained
%by surgery on a link $L$ with $m$ components in $S^3$. Suppose
%the coloring of $\Omega$ is $c_\Omega$. If $L_1$ and $L_2$
%are two framed links equipped with colorings $c_1$ and $c_2$
%respectively then denote by
%$c_1\sqcup c_2$ the corresponding 
%coloring of $L_1\sqcup L_2$. Let
%$$
%F^{P\fg}_{(L, \Omega)}(\xi) :=
%\sum_{c \in \bar{C}_k^m \cap \cC_L^R} 
%Q^\fg_{L\sqcup\Omega} (c\sqcup c_\Omega) 
%|_{q = \xi},
%$$
%where $k=r-h$.
%Then according to \cite{tu1}
%\begin{eqnarray}
%<M> & = & (F_{U_-}^{P\fg}(\xi))^\sigma \eta^{-(\sigma+m+1)} 
%F^{P\fg}_{(L, \Omega)}(\xi) \kappa^w
%\nonumber \\
%& = & \eta^{-(m+1)} F^{P\fg}_{(L, \Omega)}(\xi)
%\kappa^{w+\sigma} \nonumber \\
%& = & \eta^{-1} \frac{F^{P\fg}_{(L, \Omega)}(\xi)}
%{(F^{P\fg}_{U_-}(\xi))^{\sigma_- 
%+ \beta_1(M)}
%(F^{P\fg}_{U_+}(\xi))^{\sigma_+}}
%\kappa^{\beta_1(M) + w}. \nonumber
%\end{eqnarray}
%It is implicitly proved in \cite{le2} that
%$\frac{F^{P\fg}_{(L, \Omega)}}{F_{U_-}^{\sigma_- + \beta_1(M)}
%F_{U_+}^{\sigma_+}}$ belongs to $\BZ [\xi]$. Here $\sigma, \sigma_-$
%and $\sigma_+$ are the signature, the number of positive eigenvalues
%and the number of negative eigenvalues of $L$.
%\end{proof}

\begin{thm}\label{th_p_mfd2}
Let $M$ be an $r$-periodic integral homology 3-sphere
where $r > d h^\vee$ is a prime number.
Let $\xi$
be a primitive $r$-th root of unity then
$$
\tau_\xi^{P\fg}(M) \equiv \xi^v \overline{\tau_\xi^{P\fg}(M)} ~\mod~ 
r ~\inn~ \BZ[\xi]
$$
for some integer $v$ if $r$ is not a factor of $|G| |W|$.
\end{thm}

\begin{proof}
Lemma~\ref{lm_almost} shows the TQFT's constructed following Turaev and 
the second author \cite{tu1, le2, chenle2} are almost
integral and $\eta <M> \in \cD$.
Then the proofs given in \cite{gkp} can be used verbatim 
here to show that
the congruence relation is true in $\cD$, i.e., 
$$
x = \tau_\xi^{P\fg}(M) - \xi^v \overline{\tau_\xi^{P\fg}(M)} = r g
$$
for some $g \in \cD$. According to \cite{le2},
$x$  is in $\BZ [\xi]$. 
So $g = x/r \in \BQ (\xi)$ which implies that $g \in \BZ [\xi]$
since $\cD$ is a ring of algebra integers. 
\end{proof}

\noindent {\em Remark.}
Although one can define quantum invariants of 3-manifolds without using
TQFT the proof of theorem~\ref{th_p_mfd2} depends on it.
It would be interesting to find a direct proof.

\section{Applications}\label{ap}

It is known that the Poincare sphere $P$, -1 surgery on the left-hand
trefoil, has periodicity 2, 3 and 5. It is also known that the Brieskorn
sphere $B = \Sigma(2, 3, 7)$, -1 surgery on the right-hand trefoil, has 
periodicity
2, 3 and 7. Actually these are the only prime periodicity
they have. According to the second author \cite{le2},
$$
\tau_P = 
\tau_\xi^{Psl_2}(P) = (1 - \xi)^{-1} \sum_{n=0}^\infty \xi^n (1 -
\xi^{n+1}) (1 - \xi^{n+2}) \cdots (1 - \xi^{2n+1})
$$
and
$$ 
\tau_B = \tau_\xi^{Psl_2}(B) = 
(1 - \xi)^{-1} \sum_{n=0}^\infty \xi^{-n(n+2)} (1 -
\xi^{n+1}) (1 - \xi^{n+2}) \cdots (1 - \xi^{2n+1}).
$$ 
Here $\xi$ is again a primitive $r$-th root of unity therefore the sums
are finite. Let's suppose that $r$ is now a prime number.
Every element $x$ in the ring $\BZ[\xi]$ can be written as
$$
x = \sum_{n=0}^{r-2}a_n(x)(1 - \xi)^n + x'(1 - \xi)^{r-1}
$$
for some $x' \in \BZ[\xi]$ and some integers $a_n(x)$. The $a_n(x)'s$ are
well-defined modulo $r$ because of the following equation \cite{o1}
$$
1 + t + \cdots + t^{r-1} = {r \choose 1} + {r \choose 2}(t-1) + \cdots +
{r \choose r}(t-1)^{r-1}.
$$

One computes easily the following integers well-defined modulo $r$.
$$
\begin{array}{ll}
a_1(\tau_P) = 6, & a_1(\xi^j \overline{\tau_P}) = -6 - j, \\
a_3(\tau_P) = 464, & a_3(\xi^{-12} \overline{\tau_P}) = -16,
\end{array}
$$
$$ 
16 + 464 = 480 = (2)^5(3)(5).
$$
Hence $P$ has periodicity 2, 3 and 5 only. Similarly,
$B$ has periodicity 2, 3 and 7 only.
$$
\begin{array}{ll}
a_1(\tau_B) = 6, & a_1(\xi^j \overline{\tau_B}) = -6 - j, \\
a_3(\tau_B) = 1064, & a_3(\xi^{-12} \overline{\tau_B}) = -280,
\end{array}
$$
$$ 
1064 + 280 = 1344 = (2)^6(3)(7).
$$

\noindent
{\em Remark.} X. Zhang \cite{z1} showed that the Poincare sphere
has only periodicity 2, 3 and 5 using the fact that
the Poincare sphere is spherical.
One can also use the recently proved Orbifold Theorem to show that
the Poincare sphere is 2, 3, 5-periodic only and $\Sigma(2, 3, 7)$
is 2, 3, 7-periodic only as remarked in \cite{gkp}.

The calculation above suggests that the first and the third 
Ohtsuki's invariants
detect very well the periodicity of homology spheres. 
Let's look at these two invariants for Brieskorn spheres 
$B_n = \Sigma(2, |n|, |2n+1|)$ where $n$ is any odd number.
It can be shown, again by the Orbifold Theorem, that 
$B_n$ is $p$-periodic if and only if $p | n(2n+1)$.
Following \cite{gkp}, $B_n$ can be obtained by -1 surgery on the ($2, n$)
torus knot. Note that our $B_n$ is the $\Sigma(2, |n|, |2n-1|)$
in \cite{gkp} with opposite orientation. 
Lawrence \cite{law1} has explicit formulas for the Ohtsuki's
invariants of $B_n$. Using her formulas one can show that if $p \ge 11$ is
a prime number then $B_n$ is $p$-periodic only if $p$ divides 
$\frac{2}{3}n^2(2n+1)(n-1)(n+1)^2$. This necessary condition is sharp 
when, for example, $n = 7$ but not in general.

\noindent Qi Chen \\
Mathematics Department \\
The Ohio State University \\
Columbus, OH 43210 \\
qichen@math.ohio-state.edu \\

\noindent Thang Le \\
School of Mathematics \\
Georgia Institute of Technology \\
Atlanta, GA 30332 \\
letu@math.gatech.edu

\end{document}